\documentstyle[12pt]{amsart}

\theoremstyle{plain} 
\newtheorem{thm}{Theorem}
\newtheorem{lem}[thm]{Lemma} 
\newtheorem{cor}[thm]{Corollary}
\newtheorem{prop}[thm]{Proposition}

\newcommand{\ep}{\varepsilon} 
\newcommand{\al}{\alpha}
\newcommand{\lam}{\lambda}
\newcommand{\om}{\omega}
\newcommand{\Om}{\Omega}
\newcommand{\Sig}{\Sigma}

\newcommand{\wl}{\mbox{Weak\/$L^1$}}

\newcommand{\ol}{\overline}

\newcommand{\N}{{\Bbb N}}
\newcommand{\R}{{\Bbb R}} 
\newcommand{\Z}{{\Bbb Z}}

\newcommand{\Q}{\mbox{$\cal Q$}}

\newcommand{\spn}{\operatorname{span}}
\newcommand{\supp}{\operatorname{supp}}

\begin{document}


\title[The normed and Banach envelopes of Weak\/$L^1$]{The normed and Banach envelopes of Weak\/$L^1$}
\author{Denny H.\ Leung} 
\address{Department of Mathematics\\ National
University of Singapore\\ Singapore 119260}
\email{matlhh@@nus.edu.sg} 






\begin{abstract} 
The space \wl\ consists of all Lebesgue measurable functions on $[0,1]$ such that 
\[ q(f) = \sup_{c > 0}c\,\lam\{t : |f(t)| > c\} \]
is finite, where $\lam$ denotes Lebesgue measure. Let $\rho$ be the gauge functional of the convex hull of the unit ball $\{f : q(f) \leq 1\}$ of the quasi-norm $q$, and let $N$ be the null space of $\rho$. The normed envelope of \wl, which we denote by $W$, is the space $(\wl/N, \rho)$. The Banach envelope of \wl, $\overline{W}$, is the completion of $W$.
We show that $\ol{W}$ is isometrically lattice isomorphic to a sublattice of $W$.  It is also shown that all rearrangement invariant Banach function spaces are isometrically lattice isomorphic to a sublattice of $W$.

\end{abstract}

\maketitle 

Let $(\Om, \Sigma, \mu)$ be a measure space. The space \wl($\mu$) consists of all (equivalence classes of almost everywhere equal) real-valued $\Sigma$-measurable functions $f$ for which the quasinorm
\[ q(f) = \sup_{c > 0}c\,\mu\{\om : |f(\om)| > c\} \]
is finite. This space arose in connection with certain interpolation results, and is of importance in harmonic analysis. If $(\Om, \Sig, \mu)$ is purely non-atomic, the maximal seminorm $\rho$ on $\wl(\mu)$ was found in \cite{CF1} and \cite{CF2} to be 
\begin{equation*} 
\rho(f) = \lim_{n\to\infty}\sup\begin{Sb} q/p > n \\ p, q > 0 \end{Sb}
                   \int_{p\leq |f|\leq q}|f|\,d\mu\,/\,\ln(q/p) .
\end{equation*}
The normed envelope of \wl($\mu$) is the normed space $W(\mu) = (\wl(\mu)/N, \rho)$, where $N$ denotes the null space of the functional $\rho$. The Banach envelope is the completion $\ol{W(\mu)}$ of $W(\mu)$. In this paper, we consider (up to measure isomorphism) only the measure space $[0,1]$ endowed with Lebesgue measure $\lam$. We denote $W(\lam)$ and $\ol{W(\lam)}$ by $W$ and $\ol{W}$ respectively.
Peck and Talagrand \cite{PT} showed that $\ol{W}$ is universal for the class of all separable Banach lattices with order continuous norm. 
Recently, Lotz and Peck \cite{LP} showed that $\ol{W}$ contains isometrically lattice isomorphic copies of certain sublattices of $\ell^\infty(L^1)$. (Here and throughout, $L^1$ means $L^1[0,1]$, up to isometric lattice isomorphism.)  From this, they deduced that every separable Banach lattice is isometrically lattice isomorphic to a sublattice of $\ol{W}$. 
In this article, we show that there is a sublattice $G$ of $\ell^\infty(\ell^\infty(L^1))/c_0(\ell^\infty(L^1))$ such that $G$, $W$, and $\ol{W}$ mutually isometrically lattice isomorphically embed in one another. It is also shown that all rearrangement invariant Banach functions spaces in the sense of \cite{LT} are isometrically lattice isomorphic to sublattices of $W$.  For further results regarding the structure of \wl($\mu$), we refer the reader to \cite{KP}. Unexplained notation and terminology on vector lattices can be found in \cite{S}. If $E$ is a Banach lattice and $I$ is an arbitrary index set, let $\ell^p(I,E)$, $1 \leq p \leq \infty$, respectively, $c_0(I,E)$, be the space consisting of all families $(x_i)_{i\in I}$ such that $x_i \in E$ for all $i$, and $(\|x_i\|)_{i\in I} \in \ell^p(I)$, respectively, $c_0(I)$. We write $\ell^p(E)$ and $c_0(E)$ for these respective spaces if the index set $I = \N$. Clearly $\ell^p(I,E)$ and $c_0(I,E)$ are Banach lattices. The cardinality of a set $A$ is denoted by $|A|$.

\section{The spaces $W$ and $\ol{W}$}

If $f$ is a real-valued function defined on a set $\Om$, let the {\em support}\/ of $f$ be the set $\supp f = \{\om \in \Om : |f(\om)| > 0\}$. Furthermore, for real numbers $p\leq q$, we write $\{p \leq f \leq q\}$ for the set $\{\om \in \Om : p \leq f(\om) \leq q\}$. 

\begin{lem} \label{1}
Let $(h_k)$ be a sequence of disjointly supported Lebesgue measurable functions on $[0,1]^2$. Suppose there exist $\delta, \gamma > 0$ and strictly positive sequences $(\alpha_k)$, $(\beta_k)$ such that
\begin{enumerate}
\item $\al_k < \beta_k < \al_{k+1}$ for all $k$,
\item $\lim_k\al_k = \lim_k\beta_k = \infty$,
\item $\ln(\al_{k+1}/\beta_k) \geq (k+1)\sum^k_{j=1}\int h_j$ for all $k$, and 
\item $\delta\al_k \leq h_k(s,t) \leq \gamma\beta_k$ for all $(s,t) \in \supp h_k$.
\end{enumerate}
If $1 \leq p < q < \infty$, $q/p > \delta\al_N$, and $h$ denotes the pointwise sum $\sum h_k$, then
\[ \int_{p\leq h\leq q}h \leq \frac{1}{N}\ln\frac{\gamma q}{\delta p} + \sup_k\int_{p\leq h_k\leq q}h_k.\]
\end{lem}

\begin{pf} 
If $[\delta\al_k,\gamma\beta_k]\cap [p,q] = \emptyset$, $\int_{p\leq h_k\leq q}h_k = 0$. So we may assume that the said intersection is non-empty for some $k$. Since $\delta\al_k \to \infty$, $[\delta\al_k,\gamma\beta_k]\cap [p,q] \neq \emptyset$ for at most finitely many $k$. Let $m$ and $n$ be the minimum and maximum of the set $\{k : [\delta\al_k,\gamma\beta_k]\cap [p,q] \neq \emptyset\}$ respectively. We consider two cases.\\

\noindent\underline{Case 1} $m = n$.\\

\noindent In this case,
\[ \int_{p\leq h\leq q}h = \int_{p\leq h_m\leq q}h_m \leq \sup_k\int_{p\leq h_k\leq q}h_k .\]
\bigskip
\noindent\underline{Case 2} $m < n$.\\

\noindent Note that $p \leq \gamma\beta_m$, and $q \geq \delta\al_n$. Therefore, 
\[ \ln\frac{\gamma q}{\delta p} \geq \ln\frac{\al_n}{\beta_{n-1}} \geq n\sum^{n-1}_{k=1}\int h_k .\]
Now $q > \delta p\al_N \geq \delta\al_N$; hence $n \geq N$. Thus
\begin{align*}
\int_{p\leq h\leq q}h &= \sum^{n-1}_{k=m}\int_{p\leq h_k\leq q}h_k + \int_{p\leq h_n\leq q}h_n\\
&\leq \sum^{n-1}_{k=1}\int h_k + \int_{p\leq h_n\leq q}h_n\\
&\leq \frac{1}{n}\ln\frac{\gamma q}{\delta p} + \sup_k\int_{p\leq h_k\leq q}h_k \\
&\leq \frac{1}{N}\ln\frac{\gamma q}{\delta p} + \sup_k\int_{p\leq h_k\leq q}h_k.
\end{align*}
\end{pf}

Write any element $g \in \ell^\infty(\ell^\infty(L^1))$ as $g = (g_{ij})$, where $g_{ij} \in L^1$ for all $i, j$, and $\sup_{i,j}\|g_{ij}\|_{L^1} < \infty$. For any double sequence of numbers $M = (M_{ij})$ such that $M_{ij} \geq 1$ for all $i, j$, let $F = F_M$ be the (non-closed) lattice ideal of $\ell^\infty(\ell^\infty(L^1))$ consisting of all $g = (g_{ij}) \in \ell^\infty(\ell^\infty(L^1))$ such that $\sup_{i,j}\|g_{ij}\|_{L^\infty}/M_{ij} < \infty$. For the next result, we follow the idea of Lotz and Peck \cite{LP} in considering the \wl\ space defined on the unit square $[0,1]^2$ endowed with Lebesgue measure. Since $[0,1]$ and $[0,1]^2$ are isomorphic measure spaces, their corresponding \wl\ spaces are isometrically lattice isomorphic; the same holds for the respective normed and Banach envelopes.

\begin{prop} \label{t1}
There exists a lattice homomorphism $T : F \to W$ of norm $\leq 1$ which vanishes on $F \cap c_0(\ell^\infty(L^1))$.
\end{prop}

\begin{pf}
Choose positive sequences $(\ep_n)$ and $(r_i)$ so that $1 \geq \ep_n$, $\lim_n\ep_n = 0$, and $r_i > 1$, $\lim_ir_i = \infty$. For each $n$, let $E_n$ be the conditional expectation operator on $L^1$ with respect to the $\sigma$-algebra generated by $\{[\frac{m-1}{2^n},\frac{m}{2^n}) : 1 \leq m \leq 2^n\}$. If $i, j, n \in \N$, let $A_{ijn}$ be a countable set which is dense in 
\[ \{f \in E_nL^1 : \|f\|_{L^1} = 1, \quad \ep_n \leq f \leq nM_{ij}\} \]
with respect to the $L^\infty$-norm. For each $f \in A_{ijn}$, let $(a_m(f))^{2^n}_{m=1}$ be the coefficients such that 
\[ f = \sum^{2^n}_{m=1}a_m(f)\chi_{[(m-1)/2^n,m/2^n)}. \]
Note that $\ep_n \leq a_m(f) \leq 2^n$ for $1 \leq m \leq 2^n$. Arrange $\cup A_{ijn}$ into a sequence $(f_k)$. For each $k$, determine $i(k), j(k), n(k)$ such that $f_k \in A_{i(k),j(k),n(k)}$. Choose a positive sequence $(b_k)$ so that if we define $\al_k = b_k/2^{n(k)}$, and $\beta_k = M_{i(k),j(k)}r_{i(k)}b_k/\ep_{n(k)}$, then $\al_k < \beta_k < \al_{k+1}$, $\lim_k{\al_k} = \infty = \lim_k{\beta_k}$, and
\[ \ln\frac{\al_{k+1}}{\beta_k} \geq (k+1)\sum^k_{l=1}\ln r_{i(l)} .\]
Let $g = (g_{ij}) \in F$, and $k \in \N$. Define a function $h_k$ on $[0,1]^2$ by
\[ h_k(s,t) = \sum^{2^{n(k)}}_{m=1}\frac{g_{i(k),j(k)}(t)}{s}\chi_{B_{km}} ,\]
where 
\[ B_{km} = \biggl\{(s,t) : \frac{a_m(f_k)}{r_{i(k)}b_k} < s < \frac{a_m(f_k)}{b_k}, \frac{m-1}{2^{n(k)}} < t < \frac{m}{2^{n(k)}} \biggr\}. \]
The map $S$ defined by $Sg = \sum h_k$, where the sum is taken pointwise, is  a linear map from $F$ into the space of Lebesgue measurable functions on $[0,1]^2$. Notice that 
\begin{align*}
\supp h_k &\subseteq \bigcup^{2^{n(k)}}_{m=1}\biggl\{(s,t) : \frac{a_m(f_k)}{r_{i(k)}b_k} < s < \frac{a_m(f_k)}{b_k}\biggr\} \\
&\subseteq \biggl\{(s,t) : \frac{\ep_{n(k)}}{r_{i(k)}b_k} < s < \frac{2^{n(k)}}{b_k}\biggr\} \\
&\subseteq \biggl\{(s,t) : \frac{1}{\beta_k} < s < \frac{1}{\al_k}\biggr\} .
\end{align*}
Hence the $h_k$'s are pairwise disjoint. As the sets $B_{km}, 1 \leq m \leq 2^{n(k)}$, are also pairwise disjoint for each $k$, it follows readily that $S$ is a lattice homomorphism. Suppose $g \in F$, $\|g\| = \sup_{i,j}\|g_{ij}\|_{L^1} \leq 1$, let us estimate the $\rho$-norm of the function $Sg$. In the first instance, let us assume additionally that there exists $\delta > 0$ such that $g_{ij}(t) \geq \delta$ for all $i, j$, and $t$. Set $\gamma = \sup_{i,j}\|g_{ij}\|_{L^\infty}/M_{ij}$. If $(s,t) \in \supp h_k$, then
\[ \frac{\delta}{s} \leq \frac{g_{i(k),j(k)}(t)}{s} = h_k(s,t) \leq \frac{\gamma M_{i(k),j(k)}}{s}, \]
and
\[ \frac{M_{i(k),j(k)}}{\beta_k} = \frac{\ep_{n(k)}}{r_{i(k)}b_k} < s < \frac{2^{n(k)}}{b_k} = \frac{1}{\al_k} .\]
Hence
\[ \delta\al_k \leq h_k(s,t) \leq \gamma\beta_k .\]
Moreover,
\begin{align}\label{eq0}
\int h_k &= \sum^{2^{n(k)}}_{m=1}\int^{\frac{m}{2^{n(k)}}}_{\frac{m-1}{2^{n(k)}}}\int^{\frac{a_m(f_k)}{b_k}}_{\frac{a_m(f_k)}{r_{i(k)}b_k}}\frac{g_{i(k),j(k)}(t)}{s}\,ds\,dt \notag\\ 
& = \sum^{2^{n(k)}}_{m=1}\int^{\frac{m}{2^{n(k)}}}_{\frac{m-1}{2^{n(k)}}}g_{i(k),j(k)}(t)\,dt\,\ln r_{i(k)} \\ 
& = \|g_{i(k),j(k)}\|_{L^1}\ln r_{i(k)} \leq \ln r_{i(k)} . \notag
\end{align}
Therefore,
\[ \ln\frac{\al_{k+1}}{\beta_k} \geq (k+1)\sum^k_{l=1}\ln r_{i(l)} \geq (k+1)\sum^k_{l=1}\int h_l.\]
By Lemma \ref{1}, if $q/p > \delta\al_N$, and $p \geq 1$, then
\[ \int_{p\leq Sg\leq q}Sg \leq \frac{1}{N}\ln\frac{\gamma q}{\delta p} + \sup_k \int_{p\leq h_k\leq q}h_k .\]
If $q/p > \delta\alpha_N$ and $0 < p <1$, then, using Lemma \ref{1} again,
\[ \int_{p\leq Sg\leq q}Sg \leq \int_{p\leq Sg\leq 1}Sg + \int_{1\leq Sg\leq q/p}Sg \leq 1 + \frac{1}{N}\ln\frac{\gamma q}{\delta p} + \sup_k \int_{1\leq h_k\leq q/p}h_k .\]
Hence
\begin{equation}\label{eq1}
\lim_{n\to\infty}\sup\begin{Sb} q/p > n \\ p, q > 0 \end{Sb}
\int_{p \leq Sg \leq q}Sg/\ln(q/p) 
\leq \lim_{n\to\infty}\sup\begin{Sb} q/p > n \\ p, q > 0 \end{Sb}
\sup_k\int_{p\leq h_k\leq q}h_k/\ln(q/p) .
\end{equation}
Now 
\begin{align*}
\int_{p\leq h_k\leq q}h_k &\leq \int^1_0\int^{g_{i(k),j(k)}(t)/p}_{g_{i(k),j(k)}(t)/q}\frac{g_{i(k),j(k)}(t)}{s}\,ds\,dt\\
& = \|g_{i(k),j(k)}\|_{L^1}\ln\frac{q}{p} \leq \ln\frac{q}{p}. 
\end{align*}
Therefore, equation (\ref{eq1}) implies that $\rho(Sg) \leq 1$.  For a general $g = (g_{ij}) \in F$, and any $\delta > 0$, let $g' = (g'_{ij})$, where $g'_{ij} = |g_{ij}| + \delta$. By the above calculation, $\rho(Sg') \leq \|g'\| = \|g\| + \delta$. Since $S$ is a lattice homomorphism, $|Sg'| \geq |Sg|$. Thus $\rho(Sg) \leq \rho(Sg') \leq \|g\| + \delta$. As $\delta > 0$ is arbitrary, we conclude that $\rho(Sg) \leq \|g\|$. In particular, applying Lemma 1 in \cite{LP}, we see that $S$ maps into \wl. It is now clear that the map $T : F \to W$ defined by $Tg = Sg + N$ is a lattice homomorphism of norm $\leq 1$.

It remains to show that $T$ vanishes on $F \cap c_0(\ell^\infty(L^1))$. By the continuity of $T$, it suffices to show that $Tg = 0$ for all $g = (g_{ij}) \in F$ such that there exists $i_0 \in \N$ with $g_{ij} = 0$ whenever $i \neq i_0$.  As above, we may assume additionally that there exists $\delta >0$ such that $g_{i_0j}(t) > \delta$ for all $j$ and $t$. If $h_k \neq 0$, then $g_{i(k),j(k)} \neq 0$; hence $i(k) = i_0$. Using (\ref{eq0}),
\[ \int_{p\leq h_k\leq q}h_k \leq \int h_k \leq \|g\|\ln r_{i(k)} = \|g\| \ln r_{i_0}. \]
By (\ref{eq1}),
\[ \rho(Sg) \leq \lim_{n\to\infty}\sup\begin{Sb} q/p > n \\ p, q > 0 \end{Sb}
\frac{\|g\|\ln r_{i_0}}{\ln(q/p)} = 0. \]
\end{pf}

Let $Q$ be the quotient map from $\ell^\infty(\ell^\infty(L^1))$ onto $\ell^\infty(\ell^\infty(L^1))/c_0(\ell^\infty(L^1))$. Since $Q$ is a lattice homomorphism, $G = QF$ is a sublattice of $\ell^\infty(\ell^\infty(L^1))/c_0(\ell^\infty(L^1))$. 

\begin{thm} \label{t3}
There exists an isometric lattice isomorphism from $QF$ into $W$.
\end{thm}

\begin{pf}
Let $T$ be the map defined in the proof of Proposition \ref{t1}. Since $T$ vanishes on $F \cap c_0(\ell^\infty(L^1))$, there exists $R : QF \to W$ such that $T = RQ_{|F}$. Now $R$ is a lattice homomorphism, since both $T$ and $Q$ are, and $\|R\| \leq \|T\| \leq 1$. We claim that $\rho(RQg) \geq \|Qg\|$ for all $g \in F$. Suppose $g = (g_{ij}) \in F$, and $\|Qg\| = 1$. We may assume that there exist sequences of natural numbers $(i(l))$, $(j(i))$ such that $(i(l))$ increases to $\infty$, and $\|g_{i(l),j(l)}\|_{L^1} = 1$ for all $l$. Recall the sequence $(f_k)$ chosen in the proof of Proposition \ref{t1}. Given $\eta > 0$, there exists a sequence $({k(l)})$ in $\N$ such that $f_{k(l)} \in \cup_nA_{i(l),j(l),n}$,
\[ \sup_{l}\||g_{i(l),j(l)}| - f_{k(l)}\|_{L^1} \leq \eta, \quad\text{and}\quad \sup_l\frac{\|f_{k(l)}\|_{L^\infty}}{M_{i(l),j(l)}} < \infty .\]
Let $\phi_{ij} = f_{k(l)}$, and $\psi_{ij} = g_{i(l),j(l)}$ if $(i,j) = (i(l),j(l))$, $l \in \N$, and $\phi_{ij} = \psi_{ij} = 0$ otherwise. Then $\phi = (\phi_{ij})$ and $\psi = (\psi_{ij})$ are both in $F$, and $\|\phi - |\psi|\| \leq \eta$. Since $\|T\| \leq 1$,
\[ \rho(T\psi) = \rho(T|\psi|) \geq \rho(T\phi) - \eta .\]
Then
\[ |g| \geq \psi \implies |Tg| \geq T\psi \implies \rho(Tg) \geq \rho(T\psi) \geq \rho(T\phi) - \eta .\] 
For a given $l$, write $f_{k(l)} = \sum^{2^n}_{m=1}a_m\chi_{[(m-1)/2^n,m/2^n)}$ for some $(a_m)^{2^n}_{m=1}$, and some $n$. Note that $i(k(l)) = i(l)$, $j(k(l)) = j(l)$, and $n(k(l)) = n$. By definition of $T$, for $1 \leq m \leq 2^n$, $(s,t) \in B_{k(l),m}$,
$ |T\phi(s,t)| = \frac{a_m}{s}$. In particular, $b_{k(l)} < |T\phi(s,t)| < r_{i(l)}b_{k(l)}$ for $(s,t) \in \cup^{2^n}_{m=1}B_{k(l),m}$.
Therefore, 
\begin{align*}
\int_{b_{k(l)}\leq|T\phi|\leq r_{i(l)}b_{k(l)}}|T\phi| &\geq 
\sum^{2^n}_{m=1}\int\!\!\int_{B_{k(l),m}}\!\frac{a_m}{s}\,ds\,dt\\
& = \sum^{2^n}_{m=1}\frac{a_m}{2^n}\ln r_{i(l)} = 
\|f_{k(l)}\|_{L^1}\ln r_{i(l)}. 
\end{align*}
Since $\lim_lr_{i(l)} = \infty$, we see that $\rho(T\phi) \geq \limsup_{l}\|f_{k(l)}\|_{L^1} \geq  1 - \eta$. As $\eta > 0$ is arbitrary, it follows immediately that $\rho(RQg) =\rho(Tg) \geq 1$.
\end{pf}

Observe that if $M = (M_{ij})$ and $M' = (M'_{ij})$ satisfy $\sup_jM_{ij} = \sup_jM'_{ij} = \infty$ for all $i$, then each of $QF_M$ and $QF_{M'}$ is isometrically lattice isomorphic to a sublattice of the other. 
For the remainder of this section, let $M_{ij} = (i+1)j/\ln(i+1)$ for all $i, j \in \N$. The next result and Theorem \ref{t3} together show that $QF = QF_M$ is a maximal sublattice of $W$.

\begin{thm}\label{WtoQF}
There is an isometric lattice isomorphism from $W$ into $QF$.
\end{thm}

\begin{pf}
Given $f \in \wl$, and $i, j \in \N$, let $g_{ij} = f\chi_{\{j\leq|f|\leq(i+1)j\}}/\ln(i+1)$. It is easy to see that $g = (g_{ij}) \in F$, and that
\begin{equation}\label{eq3}
\|Qg\| = \limsup_{i\to\infty}\sup_j\|g_{ij}\|_{L^1} = \rho(f) .
\end{equation}
Consider the mapping $L : \wl \to QF$ defined by $Lf = Qg$. It follows from the proof of the Key Lemma 2.3 in \cite{KP} that $L$ is linear. Now (\ref{eq3}) tells us that the map $\tilde{L} : W \to QF$, $\tilde{L}(f+N) = Lf$, is well defined and a linear isometry. Also,
\[ \tilde{L}(|f+N|) = L|f| = Q|g| = |Qg| = |Lf| = |\tilde{L}(f+N)|.\]
Hence $\tilde{L}$ is the isometric lattice isomorphism sought. 
\end{pf}

\begin{thm}\label{envelopes}
There exists an isometric lattice isomorphism from $\ol{W}$ into $W$.
\end{thm}

\begin{pf}
It is easily verified that the set 
\[ D = \{Qg : g = (g_{ij}) \in \ell^\infty(\ell^\infty(L^1)),\quad \|g_{ij}\|_{L^\infty} \leq M_{ij}\quad\text{for all $i, j$} \} \]
is closed in $\ell^\infty(\ell^\infty(L^1))/c_0(\ell^\infty(L^1))$.
Let $\tilde{L} : W \to QF$ be the isometric lattice isomorphism given in Theorem \ref{WtoQF}. By definition of $\tilde{L}$, $\tilde{L}(W) \subseteq D$. Now there is a unique continuous linear extension $L^\# : \ol{W} \to \ell^\infty(\ell^\infty(L^1))/c_0(\ell^\infty(L^1))$ of $\tilde{L}$. Since $\tilde{L}(W) \subseteq D$, and $D$ is closed, $L^\#(\ol{W}) \subseteq D \subseteq QF$. Obviously, $L^\#$ is an isometric lattice isomorphism. 
Let $R : QF \to W$ be the isometric lattice isomorphism constructed in Theorem \ref{t3}, then $RL^\#$ is an isometric lattice isomorphism from $\ol{W}$ into $W$.
\end{pf} 

\section{Rearrangement invariant spaces}

In this section, we show that if $E$ is a rearrangement invariant space in the sense of \cite[\S 2a]{LT}, then $E$ is isometrically lattice isomorphic to a sublattice of $W$. The result is inspired by Theorem 4 in \cite{LP}, where it was shown that the Weak\/$L^p$ spaces defined on separable measure spaces are isometrically lattice isomorphic to sublattices of $\ol{W}$.
We provide the proof only for the rearrangement invariant spaces defined on $[0,\infty)$. The proofs for the measure spaces $[0,1]$ and $\N$ can be obtained by making some obvious adjustments. Recall that if $E$ is a rearrangement invariant space (or, more generally, a K\"{o}the function space \cite[Definition 1.b.17]{LT}), every measurable function $h$ such that $hf$ is integrable for all $f \in E$ defines a bounded linear functional $x'_h$ on $E$ by $x'_h(f) = \int fh$. Such functionals are called {\em integrals}. Every simple function generates an integral on $E$.

\begin{lem}\label{norm}
Let $E$ be a rearrangement invariant space on $[0,\infty)$, there exists a sequence of simple functions $(h_i)$ such that $\|x'_{h_i}\| \leq 1$ for all $n$, and $\|f\| = \limsup_{i\to\infty}|\int fh_i|$ for all $f \in E$.
\end{lem}

\begin{pf}
Let ${\cal F}$ be the collection of all simple functions of the form $h = \sum^k_{j=1}a_j\chi_{[c_{j-1},c_j)}$, where $k \in \N$, $(a_j)^k_{j=1}$, $(c_j)^k_{j=0}$ are rational numbers, and $0 = c_0 < c_1 < \dots < c_k$. Define ${\cal F}_1$ to be the subset
$\{h \in {\cal F} : \|x'_h\| \leq 1\}$.
We claim that for any $f \in E$, and any $\ep > 0$, there exists $h \in {\cal F}_1$ such that $|\int fh| > \|f\| - \ep$. Let $f \in E$ and $\ep > 0$ be given. By definition of rearrangement invariant spaces, there exists an integral $x'_g \in E'$ such that $\|x'_g\| \leq 1$, and  
$|x'_g(f)| = |\int fg| > \|f\| - \ep/2$. Let $(g_n)$ be a sequence of simple functions which converges to $g$ pointwise, and such that $|g_n| \leq |g|$ for all $n$.
By the Lebesgue Dominated Convergence Theorem, $\lim_n\int fg_n =  \int fg$. 
We may thus assume additionally that $g$ is a simple function. 
It is easy to see that there exists $h \in {\cal F}$ such that $|\int fh| \geq  |\int fg| - \ep/2 > \|f\| - \ep$, and that $h^* \leq g^*$, where $h^*$ and $g^*$ are the decreasing rearrangements of $|h|$ and $|g|$ respectively. Thus $\|x'_h\| \leq \|x'_g\| \leq 1$. Therefore, $h \in {\cal F}_1$, as desired.

Since ${\cal F}_1$ is countable, we can arrange for a sequence $(h_i)$ so that each element of ${\cal F}_1$ appears infinitely many times in the sequence. Clearly the sequence $(h_i)$ fulfills the conditions of the lemma.
\end{pf}

\begin{thm}\label{ri}
Every rearrangement invariant space $E$ on $[0,\infty)$ is isometrically lattice isomorphic to a sublattice of $W$.
\end{thm}

\begin{pf}
We will show that $E$ is isometrically lattice isomorphic to a sublattice of $\ol{QF_M}$ for some suitably chosen double sequence $M = (M_{ij})$. Then, by Theorem \ref{t3}, $E$ is isometrically lattice isomorphic to a sublattice of $\ol{W}$, which in turn is isometrically lattice isomorphic to a sublattice of $W$ by Theorem \ref{envelopes}. 

Let $(h_i)$ be the sequence given by Lemma \ref{norm}. Since $h_i$ is a simple function, there exists $0 < a_i < \infty$ such that $\supp h_i \subseteq [0, a_i]$. For $f \in E$, $i \in \N$, and $t \in [0,1]$, define $f_{i1}(t) = a_if(a_it)|h_i(a_it)|$. Also let $f_{ij} = 0$ for all $i \in \N$ and all $j > 1$. Clearly 
\begin{equation}\label{eq4}
\|f_{i1}\|_{L^1} = \int^{a_i}_0|f(u)h_i(u)|\,du = \int^\infty_0|f(u)h_i(u)|\,du .
\end{equation}
Thus $\|f_{i1}\|_{L^1} \leq \|f\|\cdot\|x'_{|h_i|}\| = \|f\|\cdot\|x'_{h_i}\| \leq \|f\|$ for all $i$. Hence $(f_{ij}) \in \ell^\infty(\ell^\infty(L^1))$. The map $T : E \to \ell^\infty(\ell^\infty(L^1))/c_0(\ell^\infty(L^1))$ defined by $Tf = Q(f_{ij})$ is easily seen to a lattice homomorphism. It follows from the preceding calculation that $\|T\| \leq 1$. On the other hand, by equation (\ref{eq4}),
\[ \limsup_i\sup_j\|f_{ij}\|_{L^1} = \limsup_i\|f_{i1}\|_{L^1} = \limsup_i\int|fh_i| \geq \limsup_i|\int fh_i| \geq \|f\|. \]
Therefore, $T$ is an isometry. 
To complete the proof, it suffices to produce a sequence $(M_i)$ such that $\lim_i\|f_{i1}\chi_{\{|f_{i1}|>M_i\}}\|_{L^1} = 0$. For then, if we define $M_{i1} = \max\{M_i,1\}$, and $M_{ij} = 1$ for $j > 1$, it is easy to check that $TE \subseteq \ol{QF_M}$, where $M = (M_{ij})$.

Let $K_i = \|h_i\|_{L^\infty}$ for all $i$. 
First note that for $f \in E$, $\|f\| \geq \int^1_0f^*(t)\,dt$; hence $c\lam\{|f|>c\} \leq \|f\|$ if $c > \|f\|$. Therefore, if $c > a_iK_i\|f\|$,
\begin{align}\label{eq5}
\lam\{|f_{i1}| > c\} &\leq \lam\bigl\{t : |f(a_it)| > \frac{c}{a_iK_i} \bigr\}\notag\\
& = \frac{1}{a_i} \lam\bigl\{|f| > \frac{c}{a_iK_i}\bigr\}\\
&\leq \frac{K_i}{c}\|f\| .\notag
\end{align}
\noindent\underline{Case 1}.  $\sup_iK_i = K < \infty$.\\

\noindent Let $(M_i)$ be any sequence such that $M_i/a_i \uparrow \infty$. Let $f \in E$. For all $i$ such that $M_i > a_iK\|f\|$, $\lam\{|f_{i1}| > M_i\} \leq K\|f\|/M_i$ by (\ref{eq5}). Hence
\begin{align*}
\|f_{i1}\chi_{\{|f_{i1}>M_i\}}\|_{L^1} &\leq \int^{K\|f\|/M_i}_0f^*_{i1}(t)\,dt \\
&\leq \int^{Ka_i\|f\|/M_i}_0f^*(t)h^*_i(t)\,dt\\
&\leq K\int^{Ka_i\|f\|/M_i}_0f^*(t)\,dt.
\end{align*}
Since $\int^1_0f^*(t)\,dt \leq \|f\| < \infty$, we obtain that $\lim_i\|f_{i1}\chi_{\{|f_{i1}|>M_i\}}\|_{L^1} = 0$.\\

\noindent\underline{Case 2}. $\sup_iK_i = \infty$.\\

\noindent For each $i$, choose $b_i > 0$ such that $h^*_i(b_i) \geq K_i/2$. Then, for all $f \in E$,
\begin{equation}\label{eq6}
\frac{K_i}{2}\int^{b_i}_0f^*(t)\,dt \leq \int^{b_i}_0f^*(t)h^*_i(t)\,dt \leq \|f\| 
\end{equation}
since $\|x'_{h^*_i}\| = \|x'_{h_i}\| \leq 1$. Let $(n_i)$ be chosen so that $\lim_iK_i/K_{n_i} = 0$. Now let $(M_i)$ be a sequence such that $(a_iK_i)^{-1}M_i > \max\{i, i/b_{n_i}\}$ for all $i$.  If $f \in E$, and $i > \|f\|$, then $\lam\{|f_{i1}| > M_i\} \leq K_i\|f\|/M_i$ by (\ref{eq5}). Therefore,
\begin{align*}
\|f_{i1}\chi_{\{|f_{i1}|>M_i\}}\|_{L^1} &\leq \int^{K_i\|f\|/M_i}_0f^*_{i1}(t)\,dt\\
&\leq K_i\int^{a_iK_i\|f\|/M_i}_0f^*(t)\,dt\\
&\leq K_i\int^{b_{n_i}}_0f^*(t)\,dt\\
&\leq \frac{2K_{i}\|f\|}{K_{n_i}} \quad \text{by (\ref{eq6})}.
\end{align*}
It follows that $\lim_i\|f_{i1}\chi_{\{|f_{i1}|>M_i\}}\|_{L^1} = 0$.
\end{pf}

Theorem \ref{ri} can be extended to certain rearrangement invariant spaces defined on non-separable measure spaces. 
Endow the two-point set $\{-1,1\}$ with the measure which assigns a mass of $1/2$ to each singleton set. For any index set $I$, denote by $\mu$ the product measure on $\{-1,1\}^I$. If $I$ is countable, $\{-1,1\}^I$ is measure isomorphic to $[0,1]$. 
For the remainder of this section, fix an index set $I$ which has the cardinality of the continuum. For each $i \in I$, let $\ep_i : \{-1,1\}^I \to \{-1,1\}$ be the projection onto the $i$-th coordinate.
If $J$ is a finite subset of $I$, and $\delta = (\delta_i)_{i\in J} \in \{-1,1\}^J$, define $\phi_{J,\delta}$ to be the function $\prod_{i\in J}\chi_{\{\ep_i = \delta_i\}}$ on $\{-1,1\}^I$. Let $\Phi_J$ be the span of the functions $\{\phi_{J,\delta} : \delta \in \{-1,1\}^J\}$.
It is not hard to see that the set $\Phi = \cup\{\Phi_J : J \subseteq I, |J| < \infty\}$ is a vector lattice (with the usual pointwise operations and order). Define $E$ by
\begin{equation}\label{E}
E = \{{\mathbf f} = (f_i)_{i\in I} : f_i \in \Phi \quad \text{for all $i$}, \quad f_i \neq 0\quad \text{for at most finitely many $i$}\}. 
\end{equation}
Similarly, let $E_J$ consist of all ${\mathbf f} = (f_i)_{i\in I} \in E$ such that $f_i \in \Phi_J$ for all $i$.
Then $E$ is a vector lattice with the coordinatewise operations and order, and $E = \cup\{E_J : J \subseteq I, |J| < \infty\}$. A norm $\|\cdot\|$ on $E$ is called a {\em lattice norm}\/ if $|{\mathbf f}| \leq |{\mathbf g}|$ implies $\|{\mathbf f}\| \leq \|{\mathbf g}\|$.
For ${\mathbf f} = (f_i) \in E$, let the {\em distribution function}\/ $d_{{\mathbf f}}$ of ${\mathbf f}$ be defined by $d_{{\mathbf f}}(t) = \sum_{i\in I}\mu\{|f_i|>t\}$, $t \geq 0$.

\begin{thm}\label{nonsep}
Let $\|\cdot\|$ be a lattice norm on $E$ which is rearrangement invariant in the sense that ${\mathbf f}$, ${\mathbf g} \in E$, $d_{{\mathbf f}} = d_{{\mathbf g}}$ implies $\|{\mathbf f}\| = \|{\mathbf g}\|$. Then $(E,\|\cdot\|)$ is isometrically lattice isomorphic to a sublattice of $W$.
\end{thm}

Of course, it follows that the completion of $E$, $\ol{E}$, is isometrically isomorphic to a sublattice of $\ol{W}$. Since $\ol{W}$ is isometrically lattice isomorphic to a sublattice of $W$ by Theorem \ref{envelopes}, the same is true for $\ol{E}$. This leads immediately to the following corollary.

\begin{cor}
If $1 \leq p < \infty$, then $\ell^p(I, L^p(\{-1,1\}^I))$ is isometrically isomorphic to a sublattice of $W$.
\end{cor}

As indicated above, $L^1$ may be identified (as a Banach lattice) with $L^1(\{-1,1\}^{\Z})$. This identification will be in force for the rest of the section. 
For each $k \in \Z$, let $r_k : \{-1,1\}^\Z \to \{-1,1\}$ be the projection onto the $k$-th coordinate. Select a bijection $\gamma : I \to \{-1,1\}^\N$. Thus, for every $i \in I$, $\gamma(i) = (\gamma(i,k))^\infty_{k=1}$, where $\gamma(i,k) = \pm 1$ for all $i \in I$, $k \in \N$. Finally, for every $i$, pick a strictly decreasing sequence of negative integers $k_i = (k_i(m))^\infty_{m=1}$ such that
\begin{itemize}
\item for each $m$, $\{k_i(m) : i \in I\}$ has only finitely many distinct values;
\item if $i \neq i'$, then $\{m : k_i(m) = k_{i'}(m)\}$ is finite.
\end{itemize}
Given a finite subset $J$ of $I$, $\delta \in \{-1,1\}^J$, $i \in I$, and $m \in \N$, define, on $\{-1,1\}^\Z$, 
\[ \psi_{J,\delta,i,m} = 2^m\,\prod^m_{k=1}\chi_{\{r_k = \gamma(i,k)\}}\cdot\prod_{j\in J}\chi_{\{r_{k_j(m)}=\delta_j\}} .\]
The mapping $T_{J,m} : E_J \to L^1$ is defined by
\[ T_{J,m}{\mathbf f} = \sum_{i\in I}\sum_{\delta\in\{-1,1\}^J}a(i,\delta)\psi_{J,\delta,i,m} \]
for all ${\mathbf f} = (f_i)_{i\in I} \in E_J$, where $f_i = \sum_{\delta\in\{-1,1\}^J}a(i,\delta)\phi_{J,\delta}$. Let us remark that the sum over $i$ is in fact a finite sum, since $f_i = 0$ for all but finitely many $i$. It is clear that $T_{J,m}$ is linear. If $I_0$ and $J$ are finite subsets of $I$, there exists $m_0 = m_0(I_0,J) \in \N$ such that 
\begin{itemize}
\item $(\gamma(i,1),\dots,\gamma(i,m_0)) \neq (\gamma(i',1),\dots,\gamma(i',m_0))$ if $i, i' \in I_0$, $i \neq i'$,
\item $k_j(m) \neq k_{j'}(m)$ whenever $j, j' \in J$, $j \neq j'$, and $m \geq m_0$.
\end{itemize}

The following lemma is easily verified by direct computation.

\begin{lem}\label{l2}
Let $I_0, J_1$, and $J_2$ be finite subsets of $I$ such that $J_1 \subseteq J_2$, and let $m \geq m_0(I_0,J_2)$.  If 
\begin{align*}
\sum_{\delta\in\{-1,1\}^{J_1}}a(i,\delta)\phi_{J_1,\delta} &= \sum_{\eta\in\{-1,1\}^{J_2}}b(i,\eta)\phi_{J_2,\eta}, \quad \text{for all $i \in I_0$, or} \\
\sum_{i\in I_0}\sum_{\delta\in\{-1,1\}^{J_1}}a(i,\delta)\psi_{J_1,\delta,i,m} &= \sum_{i\in I_0}\sum_{\eta\in\{-1,1\}^{J_2}}b(i,\eta)\psi_{J_2,\eta,i,m} ,
\end{align*}
then for all $\eta \in \{-1,1\}^{J_2}$, and all $i \in I_0$, 
$b(i,\eta) = a(i,\delta)$, where $\delta = \eta_{|J}$.
\end{lem}

An obvious consequence of the lemma is the following proposition.

\begin{prop}\label{TJm}
Let $I_0, J_1$, and $J_2$ be finite subsets of $I$ such that $J_1 \subseteq J_2$, and let $m \geq m_0(I_0,J_2)$. If ${\mathbf f} = (f_i)_{i\in I} \in E_{J_1}$, and $f_i = 0$ for all $i \notin I_0$, then $T_{J_1,m}{\mathbf f} = T_{J_2,m}{\mathbf f}$.
\end{prop}

For each ${\mathbf f} \in E$, choose a finite subset $J({\mathbf f})$ of $I$ such that ${\mathbf f} \in E_{J({\mathbf f})}$.  Given a double sequence $(h_{mn})$ of non-negative measurable functions on $\{-1,1\}^\Z$ such that $\sup_{mn}\|T_{J({\mathbf f}),m}{\mathbf f}\cdot h_{mn}\|_{L^1} < \infty$ for all ${\mathbf f} \in E$, consider the (non-linear) mapping $T : E \to \ell^\infty(\ell^\infty(L^1))$ defined by $T{\mathbf f} = (T_{J({\mathbf f}),m}{\mathbf f}\cdot h_{mn})_{mn}$.

\begin{prop}\label{qt}
Let ${Q} : \ell^\infty(\ell^\infty(L^1)) \to \ell^\infty(\ell^\infty(L^1))/c_0(\ell^\infty(L^1))$ be the quotient map. Then ${Q}T$ is a linear lattice homomorphism.
\end{prop}

\begin{pf}
Let ${\mathbf f} = (f_i)_{i\in I}$, ${\mathbf g} = (g_i)_{i\in I} \in E$, and let $\alpha \in \R$. Choose a finite subset $I_0$ of $I$ such that $f_i = 0 = g_i$ if $i \notin I_0$. Define $J = J({\mathbf f}) \cup J({\mathbf g}) \cup J(\alpha{\mathbf f}+{\mathbf g})$. If $m \geq m_0(I_0,J)$, then, for all $n \in \N$,
\begin{xalignat*}{2}
T_{J(\alpha{\mathbf f}+{\mathbf g}),m}(\alpha{\mathbf f}+{\mathbf g})\cdot h_{mn} &= T_{J,m}(\alpha{\mathbf f}+{\mathbf g})\cdot h_{mn} && \text{by Proposition \ref{TJm}}\\
&= \alpha T_{J,m}{\mathbf f}\cdot h_{mn} + T_{J,m}{\mathbf g}\cdot h_{mn} &&  \text{by linearity of $T_{J,m}$}\\
&= \alpha T_{J({\mathbf f}),m}{\mathbf f}\cdot h_{mn} + T_{J({\mathbf g}),m}{\mathbf g}\cdot h_{mn} && \text{by Proposition \ref{TJm}}.
\end{xalignat*}
Hence ${Q}T$ is linear. Now let $J' = J({\mathbf f}) \cup J(|{\mathbf f}|)$. Note that the functions $\{\psi_{J',\eta,i,m} : i \in I_0, \eta \in \{-1,1\}^{J'}\}$ are pairwise disjoint if $m \geq m_0(I_0,J')$. Thus $T_{J',m}|{\mathbf f}| = |T_{J',m}{\mathbf f}|$ for all $m \geq m_0(I_0,J')$. For all such $m$, and all $n \in \N$, it follows from Proposition \ref{TJm} that
\[ |T_{J({\mathbf f}),m}{\mathbf f}\cdot h_{mn}| = |T_{J({\mathbf f}),m}{\mathbf f}|\cdot h_{mn} = |T_{J',m}{\mathbf f}|\cdot h_{mn} = T_{J',m}|{\mathbf f}|\cdot h_{mn} = T_{J(|{\mathbf f}|),m}|{\mathbf f}|\cdot h_{mn} . \]
Therefore, $|{Q}T{\mathbf f}| = {Q}T|{\mathbf f}|$, as required.
\end{pf}

Given $m \in \N$, the set $K_m = \{k_i(m) : i \in I\}$ is a finite subset of negative integers. Let $K'_m = \{1, 2, \dots, m\} \cup K_m$. If $\eta = (\eta_k) \in \{-1,1\}^{K'_m}$, let $\zeta_{\eta,m}$ be the function $\prod_{k\in K'_m}\!\chi_{\{r_k=\eta_k\}}$ defined on $\{-1,1\}^{\Z}$. 
Associate with each real sequence $c = (c_\eta)_{\eta\in\{-1,1\}^{K'_m}}$ a function $h_c = \sum_{\eta\in\{-1,1\}^{K'_m}}c_\eta\zeta_{\eta,m}$. Also, for each $m$, choose subsets $I_m$ and $J_m$ of $I$ such that $|I_m| = 2^m$, and $|J_m| = |K_m|$. There exists a bijection $\pi_m : I_m \times \{-1,1\}^{J_m} \to \{-1,1\}^{K'_m}$. Given $c = (c_\eta)_{\eta\in\{-1,1\}^{K'_m}}$, define ${\mathbf h}_c = (h_{i,c})_{i\in I}$ by $h_{i,c} = \sum_{\tau\in\{-1,1\}^{J_m}}c_{\pi_m(i,\tau)}\phi_{J_m,\tau}$ for $i \in I_m$, and $h_{i,c} = 0$ otherwise.

\begin{lem}\label{c}
Let ${\mathbf f} = (f_i)_{i\in I} \in E$, and let $I_0$ be a finite subset of $I$ such that $f_i = 0$ if $i \notin I_0$. If $m \geq m_0(I_0,J({\mathbf f}))$, and $c = (c_\eta)_{\eta\in\{-1,1\}^{K'_m}}$, then there exists $\tilde{\mathbf h} = (\tilde{h}_i)_{i\in I}$, such that $d_{\tilde{\mathbf h}} = d_{{\mathbf h}_c}$, and $\|T_{J({\mathbf f}),m}{\mathbf f}\cdot h_c\|_{L^1} = \sum_{i\in I}\int|f_i\tilde{h}_i|$.
\end{lem}

\begin{pf}
Write $f_i = \sum_{\delta\in\{-1,1\}^{J({\mathbf f})}}a(i,\delta)\phi_{J({\mathbf f}),\delta}$ for all $i \in I_0$. There exist pairwise disjoint subsets $\{C_{i,\delta} : i\in I_0, \delta \in \{-1,1\}^{J({\mathbf f})}\}$ of $\{-1,1\}^{K'_m}$, each of cardinality $2^{|K_m|-|J({\mathbf f})|}$, such that $\psi_{J({\mathbf f}),\delta,i,m} = 2^m\sum_{\eta\in C_{i,\delta}}\zeta_{\eta,m}$. Then
\[ \|T_{J({\mathbf f}),m}{\mathbf f}\cdot h_c\|_{L^1} = \sum_{i\in I_0}\sum_{\delta\in\{-1,1\}^{J({\mathbf f})}}\sum_{\eta\in C_{i,\delta}}\frac{|a(i,\delta)c_\eta|}{2^{|K_m|}}. \]
Since $m \geq m_0(I_0,J({\mathbf f}))$, $|I_0| \leq 2^m$, and $|J({\mathbf f})| \leq |K_m|$. Choose subsets $I_1$ and $J'_m$ of $I$ such that $I_0 \cap I_1 = \emptyset$, $|I_0 \cup I_1| = 2^m$, $J({\mathbf f}) \subseteq J'_m$, and $|J'_m| = |K_m| = |J_m|$. For $i \in I_0$, $\delta \in \{-1,1\}^{J({\mathbf f})}$, there exists a bijection $\nu_{i,\delta} : C_{i,\delta} \to \{\tau\in\{-1,1\}^{J'_m} : \tau_{|J({\mathbf f})} = \delta\}$. Define $\tilde{h}_i = \sum_{\delta\in\{-1,1\}^{J({\mathbf f})}}\sum_{\eta\in C_{i,\delta}}c_\eta\phi_{J'_m,\nu_{i,\delta}(\eta)}$ for $i \in I_0$. Finally, there is a bijection $\beta : I_1 \times \{-1,1\}^{J'_m} \to \{-1,1\}^{K'_m} \backslash \cup\{C_{i,\delta} : i\in I_0, \delta\in\{-1,1\}^{J({\mathbf f})}\}$. Define $\tilde{h_i} = \sum_{\tau\in\{-1,1\}^{J'_m}}c_{\beta(i,\tau)}\phi_{J'_m,\tau}$ for $i \in I_1$. Then let $\tilde{h}_i = 0$ if $i \notin I_0 \cup I_1$. It is straightforward to check that $\tilde{\mathbf h} = (\tilde{h}_i)_{i\in I}$ fulfills the requirements of the lemma.
\end{pf}

For all $m \in \N$, let $B_m$ be the collection of all non-negative rational sequences $c = (c_\eta)_{\eta\in\{-1,1\}^{K'_m}}$ such that $\sum_{i\in I}\int|f_ih_{i,c}| \leq \|{\mathbf f}\|$ for all ${\mathbf f} = (f_i)_{i\in I} \in E$. Let us note that if $c \in B_m$, and $\tilde{\mathbf h} = (\tilde{h}_i)_{i\in I}$, $d_{\tilde{\mathbf h}} = d_{{\mathbf h}_c}$, then, due to the rearrangement invariance of the norm on $E$, $\sum_{i\in I}\int|f_i\tilde{h}_i| \leq \|{\mathbf f}\|$ for all ${\mathbf f} \in E$.

\begin{prop}\label{norming}
Let ${\mathbf f} = (f_i)_{i\in I} \in E$, and let $I_0$ be a finite subset of $I$ such that $f_i = 0$ for all $i \notin I_0$. For all $m \geq m_0(I_0,J({\mathbf f}))$,
\[ \sup_{c\in B_m}\|T_{J({\mathbf f}),m}{\mathbf f}\cdot h_c\|_{L^1} = \|{\mathbf f}\| .\]
\end{prop}

\begin{pf}
By Lemma \ref{c}, for any $c \in B_m$, there exists $\tilde{\mathbf h} = (\tilde{h}_i)_{i\in I}$ such that $d_{\tilde{\mathbf h}} = d_{{\mathbf h}_c}$, and 
$\|T_{J({\mathbf f}),m}{\mathbf f}\cdot h_c\|_{L^1} = \sum_{i\in I}\int|f_i\tilde{h}_i|$. The last sum is $\leq \|{\mathbf f}\|$ by the remark preceding the proposition. Conversely, for any $\ep > 0$, there exists $x' \in E'$, $\|x'\| \leq 1$ such that $|x'({\mathbf f})| > (1 - \ep)\|{\mathbf f}\|$. For $i_0 \in I_0$, and $\delta \in \{-1,1\}^{J({\mathbf f})}$, let ${\mathbf x}_{i_0,\delta} = (x_i) \in E$, where $x_i = \phi_{J({\mathbf f}),\delta}$ if $i = i_0$, and $x_i = 0$ otherwise. Set $b(i,\delta) = 2^{J({\mathbf f})}x'({\mathbf x}_{i,\delta})$ for $i \in I_0$, $\delta \in \{-1,1\}^{J({\mathbf f})}$. Write $f_i = \sum_{\delta\in\{-1,1\}^{J({\mathbf f})}}a(i,\delta)\phi_{J({\mathbf f}),\delta}$ for $i \in I_0$. Then
\[ (1 - \ep)\|{\mathbf f}\| < |x'({\mathbf f})| \leq \sum_{i\in I_0}\sum_{\delta\in\{-1,1\}^{J({\mathbf f})}}\frac{|a(i,\delta)b(i,\delta)|}{2^{J({\mathbf f})}}. \]
Hence, there exist non-negative rational numbers $c(i,\delta)$ such that $c(i,\delta) \leq |b(i,\delta)|$, and 
\[ (1 - \ep)\|{\mathbf f}\| < \sum_{i\in I_0}\sum_{\delta\in\{-1,1\}^{J({\mathbf f})}}\frac{|a(i,\delta)|c(i,\delta)}{2^{J({\mathbf f})}}. \]
Define ${\mathbf g} = (g_i)_{i\in I}$ by $g_i = \sum_{\delta\in\{-1,1\}^{J({\mathbf f})}}c(i,\delta)\phi_{J({\mathbf f}),\delta}$ for $i \in I_0$, $g_i = 0$ otherwise. If ${\mathbf p} = (p_i)_{i\in I} \in E$, define $P_{J({\mathbf f})}{\mathbf p} = (q_i)_{i\in I}$, $q_i = \sum_{\delta\in\{-1,1\}^{J({\mathbf f})}}2^{J({\mathbf f})}\int p_i\phi_{J({\mathbf f}),\delta}\cdot\phi_{J({\mathbf f}),\delta}$. By a standard argument, using the rearrangement invariance of the norm on $E$, we see that $\|P_{J({\mathbf f})}{\mathbf p}\| \leq \|{\mathbf p}\|$. Hence
\[ \sum_{i\in I_0}\int|p_ig_i| \leq |x'|(P_{J({\mathbf f})}|{\mathbf p}|) \leq \|{\mathbf p}\|. \]
From the proof of Lemma \ref{c}, there are pairwise disjoint subsets $\{C_{i,\delta} : i\in I_0, \delta \in \{-1,1\}^{J({\mathbf f})}\}$ of $\{-1,1\}^{K'_m}$, each of cardinality $2^{|K_m|-|J({\mathbf f})|}$, such that if we let $c_\eta = c(i,\delta)$ for all $\eta \in C_{i,\delta}$, $i\in I_0, \delta \in \{-1,1\}^{J({\mathbf f})}$, and $c_\eta = 0$ otherwise, then
for $c = (c_\eta)_{\eta\in\{-1,1\}^{K'_m}}$, 
\[ \|T_{J({\mathbf f}),m}{\mathbf f}\cdot h_c\|_{L^1} =  \sum_{i\in I_0}\sum_{\delta\in\{-1,1\}^{J({\mathbf f})}}\frac{|a(i,\delta)|c(i,\delta)}{2^{|J({\mathbf f})|}} > (1 -\ep)\|{\mathbf f}\|. \]
Note that $d_{\tilde{\mathbf h}_c} = d_{\mathbf g}$. Hence $\sum_{i\in I_0}\int|p_ih_{i,c}| \leq \|{\mathbf p}\|$ for all ${\mathbf p} = (p_i)_{i\in I} \in E$. Thus $c \in B_m$. Since $\ep > 0$ is arbitrary, we obtain the reverse inequality
\[ \sup_{c\in B_m}\|T_{J({\mathbf f}),m}{\mathbf f}\cdot h_c\|_{L^1} \geq  \|{\mathbf f}\| .\]
This completes the proof the proposition.
\end{pf}

We are now ready to prove Theorem \ref{nonsep}. For each $m \in \N$, $B_m$ is countable. Hence we can list the functions $\{h_c : c \in B_m\}$ as a sequence $(h_{mn})^\infty_{n=1}$. Define the map  $T : E \to \ell^\infty(\ell^\infty(L^1))$ by $T{\mathbf f} = (T_{J({\mathbf f}),m}{\mathbf f}\cdot h_{mn})_{mn}$. By Proposition \ref{qt}, ${Q}T$ is a lattice homomorphism, where ${Q} : \ell^\infty(\ell^\infty(L^1)) \to \ell^\infty(\ell^\infty(L^1))/c_0(\ell^\infty(L^1))$ is the quotient map. It follows from Proposition \ref{norming} that $QT$ is an (into) isometry. Finally, note that in the notation of Lemma \ref{c} and Proposition \ref{norming}, $T_{J({\mathbf f}),m}{\mathbf f}\cdot h_c \in \spn\{\zeta_{\eta,m} : \eta \in \{-1,1\}^{K'_m}\}$ for all $m \geq m_0(I_0,J({\mathbf f}))$. Hence $\|T_{J({\mathbf f}),m}{\mathbf f}\cdot h_c\|_{L^\infty} \leq 2^{|K'_m|}\|T_{J({\mathbf f}),m}{\mathbf f}\cdot h_c\|_{L^1}$. Thus $QT{\mathbf f} \in QF_M$, where $M = (M_{mn})$, $M_{mn} = 2^{|K'_m|}$ for all $m$ and $n$. An appeal to Theorem \ref{t3} yields the desired result.

\section{Order isometry}

Following \cite{LP}, we say that a linear operator $T$ from a Banach lattice $E$ into a Banach lattice $F$ is an {\em order isometry}\/ if $Tx \geq 0$ if and only if $x \geq 0$, and $\|Tx\| = \|x\|$ for all $x \in E$. In \cite{LP}, it is shown that if $E$ is a separable Banach lattice, and $E'$ has a weak order unit, then $E'$ is order isometric to a closed subspace of $\ol{W}$. Here, we show that the assumption that $E'$ has a weak order unit can be removed.

Let $\Gamma = \{-1,1\}^\N$. If $m \in \N$, and $\phi \in \Phi_m = \{-1,1\}^m$, let $\Gamma_\phi = \{\gamma \in \Gamma : \gamma_{|\{1,\dots,m\}} = \phi\}$. 

\begin{prop}\label{isom}
There is an order isometry from $\ell^\infty(\ell^1(\Gamma))$ onto a closed subspace of  $(\oplus\ell^1(\Phi_m))_{\ell^\infty}/(\oplus\ell^1(\Phi_m))_{c_0}$.
\end{prop}

\begin{pf}
Partition $\N$ into a sequence of infinite subsets $(L_n)^\infty_{n=1}$. If $a \in \ell^\infty(\ell^1(\Gamma))$, write $a = (a^n_\gamma)$, so that $\|a\| = \sup_n\sum_{\gamma\in\Gamma}|a^n_\gamma| < \infty$. Given $m \in \N$, and $\phi \in \Phi_m$, define $b_\phi = \sum_{\gamma\in\Gamma_\phi}a^n_\gamma$, where $n$ is the unique integer such that $m \in L_n$. Define the map $U : \ell^\infty(\ell^1(\Gamma)) \to (\oplus\ell^1(\Phi_m))_{\ell^\infty}$ by $Ta = b$, where $b = ((b_\phi)_{\phi\in\Phi_1}, (b_\phi)_{\phi\in\Phi_2}, \dots)$. Clearly $T$ is a positive linear operator. Note that if $m \in L_n$, then
\[ \sum_{\phi\in\Phi_m}|b_\phi| \leq \sum_{\phi\in\Phi_m}\sum_{\gamma\in\Gamma_\phi}|a^n_\gamma| = \sum_{\gamma\in\Gamma}|a^n_\gamma| \leq \|a\| . \]
Hence $\|T\| \leq 1$. Let ${\cal Q} : (\oplus\ell^1(\Phi_m))_{\ell^\infty} \to (\oplus\ell^1(\Phi_m))_{\ell^\infty}/(\oplus\ell^1(\Phi_m))_{c_0}$ be the quotient map. Then ${\cal Q}T$ is positive, and $\|{\cal Q}T\| \leq 1$. We claim that ${\cal Q}T$ is an order isometry. 

If ${\cal Q}Ta = {\cal Q}b \geq 0$, then $\lim_{m\to\infty}\sum\{b_\phi:\phi\in\Phi_m,b_\phi\leq 0\} = 0$. If $a \not\geq 0$, then there exists $n_0$ and $\gamma_0$ such that $a^{n_0}_{\gamma_0} < 0$. List the elements of $L_{n_0}$ is ascending order : $L_{n_0} = \{m_1 < m_2 < \dots\}$. For all $r \in \N$, let $\phi_r = {\gamma_0}_{|\{1,\dots,m_r\}}$. Then
\[ \lim_{r\to\infty}b_{\phi_r} = \lim_{r\to\infty}\sum_{\gamma\in\Gamma_{\phi_r}}a^{n_0}_\gamma = a^{n_0}_{\gamma_0} < 0 .\]
Thus
\[ \lim_{r\to\infty}\sum\begin{Sb}\phi\in\Phi_{m_r} \\ b_\phi\leq 0\end{Sb}b_\phi \leq a^{n_0}_{\gamma_0} < 0 ,\]
a contradiction. Therefore, $a \geq 0$.

Now, assume $\|a\| > 1$. Then there exists $n$ such that $\sum_{\gamma\in\Gamma}|a^n_\gamma| > 1$. Given $\ep > 0$, choose a finite subset $\Gamma_1$ of $\Gamma$ such that 
\[ \sum_{\gamma\in\Gamma_1}|a^n_\gamma| > 1 \quad \text{and} \quad \sum_{\gamma\notin\Gamma_1}|a^n_\gamma| < \ep .\]
Choose $m \in L_n$ so that if we define $\phi_\gamma = \gamma_{|\{1,\dots,m\}}$ for all $\gamma \in \Gamma$, then $\phi_\gamma \neq 
\phi_{\gamma'}$ for all $\gamma, \gamma' \in \Gamma_1$, $\gamma \neq \gamma'$. For $\tilde{\gamma}\in \Gamma_1$,
\[ |b_{\phi_{\tilde{\gamma}}}| = |\sum_{\gamma\in\Gamma_{\phi_{\tilde{\gamma}}}}a^n_\gamma| \geq |a^n_{\tilde{\gamma}}| - \sum\begin{Sb}\gamma\notin\Gamma_1 \\ \gamma\in\Gamma_{\phi_{\tilde{\gamma}}}\end{Sb}|a^n_\gamma| .\]
Therefore,
\[ \sum_{\gamma\in\Gamma_1}|b_{\phi_\gamma}| \geq \sum_{\gamma\in\Gamma_1}|a^n_{\gamma}| - \sum_{\gamma\notin\Gamma_1}|a^n_\gamma| > 1 - \ep .\]
Since $\ep > 0$ is arbitrary, $\|{\cal Q}Ta\| \geq \|a\|$. Since $\|{\cal Q}T\| \leq 1$ as well, we conclude that ${\cal Q}T$ is an isometry.
\end{pf}

\begin{lem}\label{L1}
Let $E$ be a separable Banach lattice. Then $E'$ is isometrically lattice isomorphic to a sublattice of $\ell^\infty(\ell^1(\Gamma,L^1))$.
\end{lem}

\begin{pf}
By the proof of Lemma 3 in \cite{LP}, for any $x \in E$, $x > 0$, there exist a compact Hausdorff space $K$, and a lattice homomorphism $S : C(K) \to E$ such that $S'$ is a lattice homomorphism, $[0,S'x']$ is weakly (and hence norm) separable, and $\|S'x'\| = |x'|(x)$ for all $x' \in E'$. Note that $E'$, and hence $S'E'$, has a dense subset of cardinality $\leq |\Gamma|$. Since $S'E'$ is a sublattice of the $AL$-space $M(K)$, has separable order intervals, and has density $\leq |\Gamma|$, it follows that $S'E'$ is isometrically lattice isomorphic to a sublattice of $\ell^1(\Gamma,L^1)$. Now let $(x_n)$ be a positive sequence in the unit ball of $E$ such that $\|x'\| = \sup_n|x'|(x_n)$ for all $x' \in E'$. For each $n$, there exists a lattice homomorphism $R_n : E' \to \ell^1(\Gamma,L^1)$ such that $\|R_nx'\| = |x'|(x_n)$ for all $x' \in E'$. Clearly, the map $R : E' \to \ell^\infty(\ell^1(\Gamma,L^1))$ defined by $Rx' = (R_nx')^\infty_{n=1}$ is an isometric lattice isomorphism.
\end{pf}

\begin{thm}
Let $E$ be a separable Banach lattice. Then $E'$ is order isometric to a closed subspace of $W$.
\end{thm}

\begin{pf}
For any $n \in \N$, let $E_n$ be the conditional expectation operator on $L^1$ with respect to the $\sigma$-algebra generated by the sets $\{[(k-1)/2^n,k/2^n) : 1 \leq k \leq 2^n\}$. Then the map $V_n : \ell^1(\Gamma,L^1) \to \ell^1(\Gamma,E_nL^1)$ defined by $V_n((f_\gamma)_{\gamma\in\Gamma}) = ((E_nf_\gamma)_{\gamma\in\Gamma})^\infty_{n=1}$ is an order isometry. Since $\ell^1(\Gamma,E_nL^1)$ is clearly isometrically lattice isomorphic to $\ell^1(\Gamma)$, it follows that $\ell^\infty(\ell^1(\Gamma,L^1))$, and hence $E'$, is order isometric to a closed subspace of $\ell^\infty(\ell^1(\Gamma))$, which in turn is order isometric to a closed subspace of $(\oplus\ell^1(\Phi_m))_{\ell^\infty}/(\oplus\ell^1(\Phi_m))_{c_0}$ by Proposition \ref{isom}. It is a simple exercise to check that the latter space is isometrically lattice isomorphic to a sublattice of $QF_M$ for a suitably chosen $M = (M_{ij})$. Finally, $QF_M$ is isometrically lattice isomorphic to a sublattice of $W$ by Theorem \ref{t3}. 
\end{pf}


\end{document}